\documentclass[12pt,draftcls,onecolumn]{IEEEtran}
\def\alexo#1{{{#1}}}
\usepackage{url}
\usepackage{color}
\usepackage{epstopdf}

\usepackage{graphicx}
\usepackage{epsfig}

\newtheorem{theorem}{Theorem}
\newtheorem{definition}[theorem]{Definition}
\newtheorem{assumption}[theorem]{Assumption}
\newtheorem{lemma}[theorem]{Lemma}

\def\E{\mathcal{E}}

\def\Q{\mathcal{Q}}
\def\A{\mathcal{A}}
\def\P{\Sigma}
\def\Ps{P_{\sigma}}
\def\e{{\bf e}}
\def\ei{{\e^{(i)}}}
\def\ej{{\e^{(j)}}}
\def\ek{{\e^{(k)}}}
\def\B{{\mathcal B}}
\def\Pm{{\mathcal P}}

\title{ \LARGE \bf
On the Nonexistence of Quadratic Lyapunov Functions for Consensus Algorithms}

\author{Alex Olshevsky, John N.
Tsitsiklis
\thanks{Laboratory for Information and Decision Systems,
Massachusetts Institute of Technology, Cambridge, MA 02139, USA;
{\tt\small alex\_o@mit.edu,
jnt@mit.edu}}%
\thanks{This research was supported by the National Science
Foundation under a Graduate Research Fellowship and grants
ECS-0312921, ECCS-0701623.} }

\begin{document}
\maketitle
\begin{abstract}
\noindent
We provide an example proving that there exists no quadratic Lyapunov function for a certain class of linear agreement/consensus algorithms, a fact that had been numerically verified in \cite{JLM03}. We also briefly discuss sufficient conditions for the existence of such a Lyapunov function. \end{abstract}

\section{Introduction}

We examine a class of algorithms that can be used by a group of
agents (e.g., UAVs, nodes of a communication network, etc.) in order
to reach consensus on a common opinion (represented by a scalar or
vector), starting from different initial opinions, and possibly in
the presence of severe restrictions on inter-agent communications.

We focus on a particular algorithm, whereby, at each time step, every agent averages its own opinion with received messages containing the current opinions of some other agents. While this algorithm is known to converge under mild conditions, convergence proofs usually rely on the ``span norm''
of the vector of opinions. In this note, we address the question of whether convergence can also be established using a quadratic Lyapunov function. Among other reasons, this question is of interest because of its potential implications on convergence time analysis.
A negative answer to this question was provided in \cite{JLM03}, where the nonexistence
of a quadratic Lyapunov function was verified numerically. In this paper, we provide an explicit example and proof of this fact.

In Section II we give some definitions and formally state the problem.
Section III contains the main result and its proof. Section IV provides some additional perspective, together with some conditions under which a quadratic Lyapunov function is guaranteed to exist.

\section{The Agreement Algorithm}

We consider a set $N=\{1,2,\ldots,n \}$ of agents embedded, at each nonnegative integer time $t$, in a directed graph $G(t)=(N,\E(t))$. We assume that $(i,i)\in \E(t)$, for all $i$ and $t$. We define $N_i(t)=\{j\mid (j,i)\in\E(t)\}$, and let $d_i(t)$ be the cardinality of $N_i(t)$.

Each agent $i$ starts with a scalar value $x_i(0)$. At each time $t$, agent $i$ receives from every agent $j\in N_i(t)$ a message with the value of $x_j(t)$,  and uses the received values to perform the update
$$ x_i(t+1)=\sum_{j=1}^n a_{ij}(t)x_j(t),
$$
where the $a_{ij}(t)$ are nonnegative coefficients that satisfy
$a_{ij}(t)=0$ if $(j,i) \notin \E(t)$, and $\sum_{j\in N_i(t)}
a_{ij}(t)=1$, so that $x_i(t+1)$ is a weighted average of the values
$x_j(t)$ held by the agents at time $t$. We define the vector
$x(t)=(x_1(t),\ldots,x_n(t))$, and note that the algorithm can be
written in the form $x(t+1)=A(t)x(t)$.

We next state some conditions under which the agreement algorithm is
guaranteed to converge.
\vspace{5pt}

\begin{assumption} There exists some $\alpha>0$ such that if $(j,i)\in\E(t)$, then $a_{ij}(t)\geq \alpha$.
\end{assumption}

\begin{assumption} (Bounded intercommunication intervals)
There is some $B$ such that for every nonnegative integer $k$, the
graph $(N,\mathcal{E}(kB)\cup \mathcal{E}(kB+1)\cup\cdots \cup
\mathcal{E}((k+1)B))$ is strongly connected.
\end{assumption}

\begin{theorem} \label{mainth}
Under Assumptions 1-2, and for every $x(0)$, the components $x_i(t)$, $i=1,\ldots,n$, converge to a common limit.
\end{theorem}

Theorem \alexo{\ref{mainth}} is presented in \cite{TBA86} and is
proved in \cite{T84} (under a slightly different version of
Assumption 2), as well as in \cite{JLM03}, for a special case to be
considered below; \alexo{see also \cite{HB05, M05} for
generalizations and extensions}. On the other hand, if the graphs
$\alexo{G}(t)$ are symmetric, namely, $(i,j)\in \mathcal{E}(t)$ if
and only if $(j,i)\in \mathcal{E}(t)$, Assumption 2 can be replaced
by the weaker requirement that the graph $(N,\cup_{s\geq t}
\mathcal{E}(t))$ is strongly connected for every $t\geq 0$; see
\cite{HB05,LW04,CMA04,M05}.

We will focus on a special case, motivated from the model of Vicsek
et al.\ \cite{VCBJCS95}, and studied in \cite{JLM03}, to be referred
to as the {\it symmetric, equal-neighbor,} model. In this model, the
graphs $G(t)$ are symmetric, and $a_{ij}(t)=1/d_i(t)$, for every
$(j,i)\in\E(t)$. Thus, each node $i$ forms an unweighted average of
the values $x_j(t)$ that it has access to (including its own).

Theorem 1 is usually proved by showing that the ``span norm''
$\max_{i} x_i(t) - \min_i x_i(t)$ is guaranteed to decrease after a certain number
of iterations. Unfortunately, this proof method usually gives an overly conservative bound on the convergence time of the algorithm.
Tighter bounds on the convergence time would have to rely on alternative Lyapunov functions, such as quadratic ones, of the form $x^T M x$, if they exist.

\alexo{Although quadratic Lyapunov functions can always be found for
linear systems, they may fail to exist when the system is allowed to
switch between a fixed number of linear modes. On the other hand,
there are classes of such switched linear systems that do admit
quadratic Lyapunov functions. See \cite{LM99} for a broad overview
of the literature on this subject. For the symmetric, equal-neighbor
model this issue was investigated in \cite{JLM03}.} The authors
write:

\begin{quote} ``...no such common Lyapunov matrix $M$ exists. While we have not been able to construct a simple analytical example which demonstrates
this, we have been able to determine, for example, that no common
quadratic Lyapunov function exists for the class of all [graphs
which have] $10$ vertices and are connected. One can verify that
this is so by using semidefinite programming...'' \end{quote}

The main contribution of this note is to provide an analytical
example that proves this fact.

\section{The Example}

Let us fix a positive integer $n$. We start by defining a class
$\alexo{\Q}$ of functions with some minimal desired properties of
quadratic Lyapunov functions. Let $\e$ be the vector in $\Re^n$ with
all components equal to 1. A square matrix is said to be {\it
stochastic} if it is nonnegative and the sum of the entries in each
row is equal to one. Let $\A \subset \Re^{n \times n}$ be the set of
stochastic matrices $A$ such that: (i) $a_{ii}>0$, for all $i$; (ii)
all positive entries on any given row of $A$ are equal; (iii)
$a_{ij} > 0$ if and only if $a_{ji}>0$; (iv) the graph associated
with the set of edges $\{ (i,j) ~|~ a_{ij}>0 \}$ is connected. These
are precisely the matrices that correspond to a single iteration of
the equal-neighbor algorithm on symmetric, connected graphs.

\begin{definition} \label{lyap-def}
\label{d:l} A function $\alexo{Q}:\Re^n\to\Re$ belongs to the class
$\alexo{\Q}$ if it is of the form $\alexo{Q}(x)=x^TMx$, where:
\begin{itemize}
\item[(a)]
The matrix $M \in \Re^{n \times n}$ is nonzero, symmetric, and
nonnegative definite.
\item[(b)] For every $A\in\A$, and $x\in\Re^n$, we have
$\alexo{Q}(Ax) \leq \alexo{Q}(x)$.
\item[(c)] We have $\alexo{Q}(\e)=0$.
\end{itemize}
\end{definition}

Note that condition (b) may be rewritten in matrix form as
\begin{equation} \label{cond:b} x^TA^TMAx \leq x^T M x, \qquad
\mbox{for all } A \in \A, \mbox{ and } x \in \Re^n.
\end{equation} The rationale behind condition (c) is as follows.
Let $S$ be the subspace spanned by the vector $\e$. Since we are
interested in convergence to the set $S$, and every element of $S$
is a fixed point of the algorithm, it is natural to require that
$\alexo{Q}(\e)=0$, or, equivalently,
$$M\e=0.$$
Of course, for a Lyapunov function to be useful, additional
properties would be desirable. For example we should require some
additional condition that guarantees that $\alexo{Q}(x(t))$
eventually decreases. However, according to Theorem~\ref{th:main},
even the minimal requirements in Definition \alexo{\ref{lyap-def}}
are sufficient to preclude the existence of a quadratic Lyapunov
function.

\begin{theorem} \label{th:main}
Suppose that $n \geq 8$. Then, the class $\alexo{\Q}$ (cf.\
Definition \ref{d:l}) is empty.
\end{theorem}

The idea of the proof is as follows. Using the fact the dynamics of
the system are essentially the same when we rename the components,
we show that if $x^TMx$ has the desired properties, so does $x^T
\alexo{Z} x$ for a matrix $\alexo{Z}$ that has certain
permutation-invariance properties. This leads us to the conclusion
that there is essentially a single candidate Lyapunov function, for
which a counterexample is easy to develop.

Recall that a permutation of $n$ elements is a bijective mapping
$\sigma:\{1,\ldots,n\}\to \{1,\ldots,n\}$. Let $\P$ be the set of
all permutations of $n$ elements. For any $\sigma\in\P$, we
define a corresponding permutation matrix $P_{\sigma}$ by letting
the $i$th component of $P_{\sigma}x$ be equal to $x_{\sigma(i)}$.
Note that $\Ps^{-1}=\Ps^T$, for all $\sigma\in \P$. Let $\Pm$ be the
set of all permutation matrices corresponding to permutations in
$\P$.
\begin{lemma} \label{sumpermutations} Let $M \in \alexo{\Q}$. Define $\alexo{Z}$ as
\[ \alexo{Z} = \sum_{P\in\Pm} P^T M P. \] Then, $\alexo{Z} \in \alexo{\Q}$.
\end{lemma}

\noindent {\bf Proof:} For every matrix $A\in\A$, and any $P\in
\Pm$, it is easily seen that $P A P^T \in\A$. This is because
the transformation $A\mapsto P A P^T$ amounts to permuting the
rows and columns of $A$, which is the same as permuting (renaming)
the nodes of the graph.

We claim that if $M\in \alexo{\Q}$ and $P\in\Pm$, then $P^T M P\in
\alexo{\Q}$. Indeed, if $M$ is nonzero, symmetric, and nonnegative
definite, so is $P^T M P$. Furthermore, since $P\e =\e$, if $M\e=0$,
then $P^T M P \e= 0$. To establish condition (b) in Definition
\ref{d:l}, let us introduce the notation $\alexo{Q}_P(x) = x^T (P^T
M P) x$. Fix a vector $x\in\Re^n$, and $A \in \A$; define $B=P A
P^T\in\A$. We have
\begin{eqnarray*} \alexo{Q}_P(Ax) & = & x^T A^T P^T M P A x \\
&=&x^T P^T P A^T P^T M P A P^T P x\\
&=&x^T P^T B^{T} M B P x\\
&\leq& x^T P^T M P x\\
& = & \alexo{Q}_P(x),
\end{eqnarray*}
where the inequality follows by applying Eq. (\ref{cond:b}), which
is satisfied by $M$, to the vector $Px$ and the matrix $B$. We
conclude that $\alexo{Q}_P \in \alexo{\Q}$.

Since the sum of matrices in $\alexo{\Q}$ remains in $\alexo{\Q}$,
it follows that $\alexo{Z}= \sum_{P\in\Pm} P^T M P$ belongs to
$\alexo{\Q}$. \rule{7pt}{7pt}

We define the ``sample variance'' $V(x)$ of the values
$x_1,\ldots,x_n$, by
$$V(x)= \sum_{i=1}^n (x_i - \bar{x})^2,$$ where
$\bar{x} = (1/n) \sum_{i=1}^n x_i$. This is a nonnegative quadratic
function  of $x$, and therefore, $V(x)=x^T C x$, for a suitable
nonnegative definite, nonzero symmetric matrix $C \in \Re^{n \times
n}$.

\begin{lemma} \label{qvariance} There exists some $\alpha>0$ such that
\[ x^T \alexo{Z} x = \alpha V(x), \qquad \mbox{for all } x\in \Re^n. \] \end{lemma}

\noindent {\bf Proof}: We observe that the matrix $\alexo{Z}$
satisfies
\begin{equation}
\alexo{R}^T \alexo{Z} \alexo{R} = \alexo{Z}, \qquad \mbox{for all }
R \in \Pm. \label{eq:p}
\end{equation}
To see this, fix $R$ and notice that the mapping $P\mapsto P R$ is a
bijection of $\Pm$ onto itself, and therefore,
$$R^T \alexo{Z} R= \sum_{P\in\Pm} (PR)^T M (P R)
=\sum_{P\in\Pm} P^T M P =\alexo{Z}.$$

We will now show that condition (\ref{eq:p}) determines $\alexo{Z}$,
up to a multiplicative factor. Let $\alexo{z}_{ij}$ be the $(i,j)$th
entry of $\alexo{Z}$. Let $\ei$ be the $i$th unit vector, so that
$\ei^T \alexo{Z} \ei = \alexo{z}_{ii}$. Let $P\in\Pm$ be a
permutation matrix that satisfies $P\ei=\ej$. Then,
$\alexo{z}_{ii}=\ei^T \alexo{Z} \ei =\ei^T P^T \alexo{Z} P \ei=\ej^T
 \alexo{Z} \ej =\alexo{z}_{jj}$. Therefore, all diagonal entries of $\alexo{Z}$ have a common
value, to be denoted by $\alexo{z}$.

Let us now fix three distinct indices $i,j,k$, and let $y=\ei+\ej$,
$\alexo{w}=\ei+\ek$.  Let $P\in\Pm$ be a permutation matrix such
that $P\ei=\ei$ and $P\ej=\ek$, so that $Py=\alexo{w}$. We have
$$
2 \alexo{z} +2 \alexo{z}_{ij} = y^T \alexo{Z} y= y^T P^T \alexo{Z} P
y = \alexo{ w^T Z w =2z + 2 z_{ik}}.$$ By repeating this argument
for different choices of $i,j,k$, it follows that all off-diagonal
entries of $\alexo{Z}$ have a common value to be denoted by $r$.
Using also the property that $\alexo{Z} \e=0$, we obtain that
$\alexo{z}+(n-1)r=0$. This shows that the matrix $\alexo{Z}$ is
uniquely determined, up to a multiplicative factor.

We now observe that permuting the components of a vector $x$ does
not change the value of $V(x)$. Therefore, $V(x)=V(Px)$ for every
$P\in\Pm$, which implies that $x^TP^T CPx=x^T C x$, for all $P\in
\Pm$ and $x\in\Re^n$. Thus, $C$ satisfies (\ref{eq:p}). Since all
matrices that satisfy (\ref{eq:p}) are scalar multiples of each
other, the desired result follows. \rule{7pt}{7pt}
\begin{center}
\begin{figure}
\hspace{2cm}
\includegraphics[width=10cm]{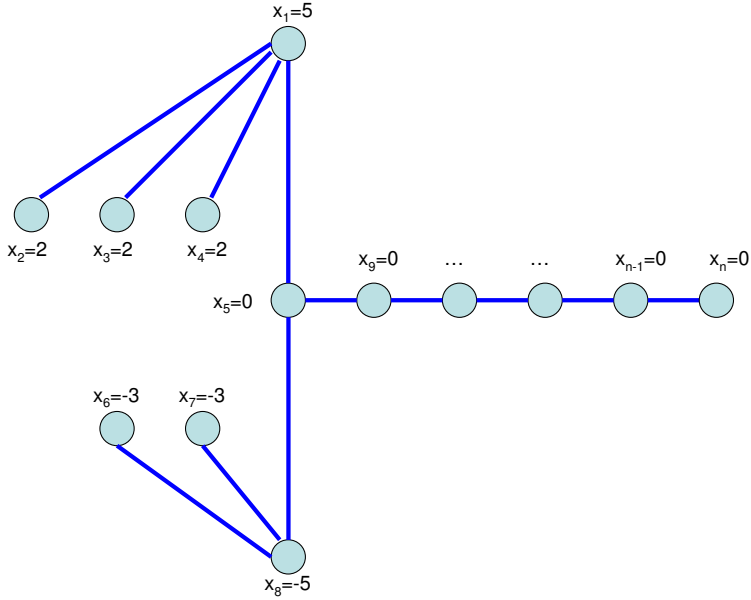}
\caption{\label{f1}
 A connected graph
on $n$ nodes showing that $V(x)$ is not a Lyapunov function when $n
\geq 8$. All arcs of the form $(i,i)$ are also assumed to be
present, but are not shown. The nodes perform an iteration of the
symmetric, equal-neighbor model according to this graph.}
\end{figure}
\end{center}

\noindent {\bf Proof of Theorem \ref{th:main}:} In view of Lemmas
\ref{sumpermutations} and \ref{qvariance}, if $\alexo{\Q}$ is
nonempty, then $V \in \alexo{\Q}$. Thus,  it suffices to show that
$V \notin \alexo{\Q}$. Suppose that $n \geq 8$, and consider the
vector $x$ with components $x_1=5$, $x_2=x_3=x_4=2$, $x_5=0$,
$x_6=x_7=-3$, $x_8=-5$, and $x_9 = \cdots = x_n = 0$. We then have
$V(x)=80$. Consider the outcome of one iteration of the symmetric,
equal-neighbor algorithm, if the graph has the form shown in Figure
\ref{f1}. After the iteration, we obtain the vector $y$ with
components  $y_1 = 11/5, y_2=y_3=y_4=7/2, y_5=0, y_6=y_7= -4, y_8 =
-11/4$, and $y_9=\cdots=y_n=0$. We have
\begin{eqnarray} V_n(y) & = & \sum_{i=1}^n (y_i - \bar{y})^2\nonumber \\
& \geq & \sum_{i=1}^8 (y_i - \bar{y})^2 \nonumber\\
& \geq & \sum_{i=1}^8 \Big(y_i - \frac{1}{8} \sum_{i=1}^8 y_i\Big)^2,
\label{eq:v8}
\end{eqnarray}
where we used that $\sum_{i=1}^k (y_i - z)^2$ is minimized
          when $z=(1/k) \sum_{i=1}^k y_i$. A
simple calculation shows that
the expression (\ref{eq:v8}) evaluates to $10246/127 \approx 80.68$,
which implies that $V(y)> V(x)$. Thus, if $n\geq8$, $V \notin
\alexo{\Q}$, and the set $\alexo{\Q}$ is empty. \rule{7pt}{7pt}

\section{Conditions for the Existence of a Quadratic Lyapunov Function}
Are there some additional conditions (e.g., restricting the matrices
$A$ to a set smaller than $\A$), under which a quadratic Lyapunov
function is guaranteed to exist? We start by showing that the answer
is positive for the case of a fixed matrix (that is, if the graph
$G(t)$ is the same for all $t$).

Let $A$ be a stochastic matrix, and suppose that there exists a
positive vector $\pi$ such that $\pi^T A = \pi^T$. Without loss of
generality, we can assume that $\pi^T \e =1$. It is known that in
this case,
\begin{equation}\label{contr}
x^T A^T D A x \leq x^T D x,\qquad \forall\ x\in\Re^n,
\end{equation}
where $D$ is a diagonal matrix, whose $i$th diagonal entry is equal
to $\pi_i$ (cf.\ Lemma 6.4 in \cite{NDP}). However, $x^T D x$ cannot
be used as a Lyapunov function because $D\e \neq 0$ (cf.\ condition
(c) in Definition \ref{d:l}). To remedy this, we argue as in
\cite{BHOT05} and define the matrix $H=I-\e \pi^T$, and consider the
choice $M=H^T D H$. Note that $M$ has rank $n-1$.

We have $H\e= (I-\e \pi^T) \e =\e -\e (\pi^T \e) =\e-\e=0$, as
desired. Furthermore,
$$HA=A -\e \pi^T A =A-\e \pi^T = A - A \e \pi^T = AH.$$ Using this property, we obtain, for every $x\in\Re^n$,
$$x^T A^T M A x = x^T A^T H^T D H A x
=(x^T H^T) A^T D A (H x) \leq x^T H^T D H x = x^T M x,$$ where the
inequality was obtained from (\ref{contr}), applied to $Hx$. This
shows that $H^T D H$ has the desired properties (a)-(c) of
Definition \ref{d:l}, provided that $\A$ is replaced with $\{A\}$.

We have just shown that every stochastic matrix (with a positive
left eigenvector associated to the eigenvalue 1) is guaranteed to
admit a quadratic Lyapunov function, in the sense of Definition
\ref{d:l}. \alexo{Moreover, our discussion implies that there are
some classes of stochastic matrices $\B$ for which the same Lyapunov
function can be used for all matrices in the class.}
\begin{itemize}
\item[(a)]
Let $\B$ be a set of stochastic matrices. Suppose that there exists
a positive vector $\pi$ such that $\pi^T \e =1$, and $\pi^T A
=\alexo{\pi^T}$ for all $A\in\B$. Then, there exists a nonzero,
symmetric, nonnegative definite matrix $M$, of rank $n-1$, such that
$M\e=0$, and $x^T A^T M Ax \leq x^T M x$, for all $x$ and $A\in\B$.
\item[(b)]
The condition in (a) above is automatically true if all the matrices
in $\B$ are doubly stochastic (recall that a matrix $A$ is doubly
stochastic if both $A$ and $A^T$ are stochastic); in that case, we
can take $\pi=\e$.
\item[(c)] The condition in (a) above holds if and only if there exists a positive vector $\pi$, such that $\pi^T A x = \pi^T x$, for all $A\in\B$ and all $x$. In words, there must be a positive linear functional of the agents' opinions which is conserved at each iteration. For the case of doubly stochastic matrices, this linear functional is any positive multiple of the sum $\sum_{i=1}^n x_i$ of the agents' values (e.g., the average of
these values).
\end{itemize}

\noindent {\bf Acknowledgments}\vspace{0.3cm}

\noindent The authors are grateful to Ali Jadbabaie for useful
discussions about this problem.


\begin{thebibliography}{99}

\bibitem{BT89} D.\ P.\ Bertsekas, and J.\ N.\ Tsitsiklis, {\it Parallel and Distributed
Computation: Numerical Methods,} Prentice Hall, 1989.

\bibitem{NDP} D.\ P.\ Bertsekas, and J.\ N.\ Tsitsiklis, {\it Neuro-Dynamic Programming,} Athena Scientific, 1996.

\bibitem{BHOT05} V. D. Blondel, J.M. Hendrickx, A. Olshevsky, J.N.
Tsitsiklis, ``Convergence in Multiagent Coordination, Consensus, and
Flocking,'' {\it  Proceedings of the Joint 44th IEEE Conference on
Decision and Control (CDC '05),} Seville, Spain, December 2005.

\bibitem{CMA04} M. Cao, A. S. Morse, B. D. O. Anderson, ``Coordination of an
Asynchronous, Multi-Agent System via Averaging,'' {\it Proceedings
of he 16th International Federation of Automatic Control World
Congress(IFAC '05),} Prague, Czech Republic, July 2005.

\bibitem{HB05} J. M. Hendrickx and V. D. Blondel, ``Convergence of different
linear and non-linear Vicsek models,'' {\it Proceedings of the 17th
International Symposium on Mathematical Theory of Networks and
Systems (MTNS 2006)}, Kyoto (Japan), July 2006.

\bibitem{JLM03} A.\  Jadbabaie, J.\  Lin, and A.\ S.\  Morse, ``Coordination of groups
of mobile autonomous agents using nearest neighbor rules,'' {\it
IEEE Transactions on Automatic Control}, Vol.\ 48, No.\  3, pp.\
988-1001, 2003.

\bibitem{LW04} S.\  Li and H.\  Wang, ``Multi-agent coordination using
nearest-neighbor rules: revisiting the Vicsek model'';
\url{http://arxiv.org/abs/cs.MA/0407021}.

\bibitem{LM99} D. Liberzon, A.S. Morse, ``Basic problems in stability and design of switched systems,''
{\it IEEE Control Systems Magazine, } 19, no. 5, pp. 59-70, 1999.

\bibitem{M05} L. Moreau, ``Consensus seeking in multi-agent systems using
dynamically changing interaction topologies,'' {\it IEEE
Transactions on Automatic Control,} vol 50, No. 2, 2005.

\bibitem{T84} J. N. Tsitsiklis,``Problems in Decentralized Decision Making
and Computation,'' Ph.D. Thesis, Department of EECS, MIT, 1984.
\url{http://hdl.handle.net/1721.1/15254}

\bibitem{TBA86} J.\ N.\ Tsitsiklis, D.\ P.\ Bertsekas, and M.\ Athans,
``Distributed Asynchronous Deterministic and Stochastic Gradient
Optimization Algorithms," {\em IEEE Transactions on Automatic
Control,} Vol.\ 31, No.\ 9,
   1986.

\bibitem{VCBJCS95} T.\ Vicsek, E.\ Czirok, E.\ Ben-Jacob, I.\ Cohen, and O.\
Shochet. ``Novel Type of Phase Transitions in a System of
Self-Driven Particles,'' {\it Physical Review Letters,} 75, 1995.

\end{thebibliography}
\end{document}